\documentclass[a4paper,draft]{amsproc}
\usepackage{amssymb}
\usepackage{amscd} 

\theoremstyle{plain}

\theoremstyle{definition}

\theoremstyle{remark}
  
 \numberwithin{equation}{section}

 \usepackage{url}
\def\C {\mathbb{C}}
\def\Q {\mathbb{Q}}

\def\Z {\mathbb{Z}}

\def\QQ {\overline{\Q}}

\renewcommand{\le}{\leqslant}
\renewcommand{\ge}{\geqslant}\renewcommand{\geq}{\geqslant}
\renewcommand{\setminus}{\smallsetminus}
\title[\MakeUppercase{The Schanuel Subset Conjecture implies Gelfond's Conjecture}]{\MakeUppercase{The Schanuel Subset Conjecture implies}\\\MakeUppercase{Gelfond's power tower conjecture}}



\author[Diego Marques]{\bfseries Diego Marques}

\address{
Department of Mathematics \\ 
University of Bras\' ilia\\ 
Brasilia-DF\\
Brazil}
\email{diego@mat.unb.br}

\author[Jonathan Sondow]{\bfseries Jonathan Sondow}
\address{
209 West 97th Street \\ 
New York\\ 
NY 10025\\
USA}
\email{jsondow@alumni.princeton.edu}

\thanks{The first author was partially supported by FEMAT-Brazil and CNPq-Brazil} 

\begin{document}

\vspace{18mm}
\setcounter{page}{1}
\thispagestyle{empty}


\begin{abstract}
As an alternative to the famous Schanuel's Conjecture (SC), we introduce the {\it Schanuel Subset Conjecture} (SSC):
\emph{Given \mbox{$\alpha_1,\dots,\alpha_n\in\C$} linearly independent over $\Q\,$, if $\{\alpha_1,\dots,\alpha_n$, $e^{\alpha_1},\dots,e^{\alpha_n}\}$ is $\overline{\Q}$-dependent on a subset $\{\beta_1,\dots,\beta_n\}$, then $\beta_1,\dots,\beta_n$ are algebraically independent}. (A set $X\subset \C$ is called \emph{$\overline{\Q}$-dependent on} $Y\subset \C$ if \mbox{$\overline{\Q}(X) \subset \overline{\Q}(Y)$}.) We discuss whether SC is equivalent to the a~priori weaker SSC. Assuming SSC, we give conditional proofs of Gelfond's Power Tower Conjecture and of two other results.

\end{abstract}

\maketitle

\section{Introduction}

In Paris at the second International Congress of Mathematicians, Hilbert \cite[pp.~15--16]{browder}, \cite{hilbert} raised the following problem as the seventh in his famous list of twenty-three.\\

\noindent{\bf Hilbert's 7th Problem (1900).}
``{\it Prove that the expression $\alpha^{\beta}$, for an algebraic base $\alpha\neq 0,1$ and an irrational algebraic exponent $\beta$, e.g., the number $2^{\sqrt{2}}$ or $e^{\pi}= i^{\,-2i}$, always represents a transcendental or at least an irrational number.}''\\

Twenty-nine years later, Gelfond \cite{gelfond29} proved the special case in which $\beta$ is an imaginary quadratic irrational. Since $e^{\pi i}=-1$ is algebraic, his result includes the following.\\

\noindent{\bf Gelfond's Theorem (1929).}
{\it The number $e^{\pi}$ is transcendental.}\\

Five years afterward, Gelfond and Schneider, independently, solved Hilbert's 7th Problem (see Baker~\cite[p.~9]{baker} and Tijdeman's essay in \cite[pp.~241--268]{browder}).\\

\noindent{\bf Gelfond-Schneider Theorem (1934).}
\textit{If $\alpha$ and $\beta$ are algebraic numbers, where $\alpha \neq 0,1$ and $\beta$ is irrational, then $\alpha^{\beta}$ is transcendental.}\\

In the same year 1934, Gelfond \cite{gelfond} announced without proof two vast extensions of the Gelfond-Schneider Theorem.

\begin{quotation}
\hspace{0.5 cm}
Par une g\'en\'eralisation de la m\'ethode qui sert pour la d\'e\-mons\-tra\-tion du th\'eor\`eme [de Gelfond-Schneider], j'ai d\'emontr\'e les th\'e\-o\-r\`e\-mes plus g\'en\'eraux suivants $\dotso$.

La d\'emonstration de ces r\'esultats et de quelques autres r\'e\-sul\-tats sur les nombres transcendants sera donn\'ee dans un autre Recueil.
\end{quotation}

However, in 2010 in reply to a question by the second author, Michel Waldschmidt wrote, ``Gelfond never published proofs, but in 1948 and later he published proofs of much weaker statements---so it is clear now that he did not have a complete proof [of the extensions].'' (For the much weaker statements, as well as proofs of them and references, see Feldman and Nesterenko \cite[pp. 260--267]{fn}.)

We formulate Gelfond's second extension as a  conjecture and quote it verbatim, including his partial italicization and his omission of the hypothesis that the $\alpha_k$ are irrational.\\

\noindent{\bf Gelfond's Second Conjecture.}
{\it ``Les nombres
\begin{equation*}
e^{\,\omega_1e^{\,\omega_2e^{\,\cdot^{\cdot^{\cdot^{\,\omega_{n-1}e^{\,\omega_n}}}}}}} \quad\emph{et}\quad \alpha_1^{\,\alpha_2^{\, \alpha_3^{\,\cdot^{\cdot^{\cdot^{\,\alpha_m}}}}}},
\end{equation*}
o\`u $\omega_1\neq0,\omega_2,\dotsc,\omega_n$ \emph{et} $\alpha_1\neq0,1,\alpha_2\neq0,1,\alpha_3\neq0,\alpha_4,\dotsc,\alpha_m$ sont des nombres alg\' ebriques, sont des nombres transcendants et entre les nombres de cette forme n'existent pas des relations alg\' ebriques, \` a coefficients rationels \emph{(non triviales)}.}''\\

Here is a corrected version of the special case with $\omega_1=\dotsb=\omega_n$ and $\alpha_1=\dotsb=\alpha_m$, which we call:\\

\noindent{\bf Gelfond's Power Tower Conjecture.}
{\it Let $\omega\neq0$ and $\alpha$ be algebraic numbers, with $\alpha$ irrational. Then when $z:=e^{\,\omega}$ and when $z:=\alpha$, the power tower of~$z$ of order $k\ge2$
\begin{equation*}
^k\!z:=\underbrace{z^{\,z^{\,z^{\,\cdot^{\cdot^{\cdot^{z}}}}}}}_{k} 
\end{equation*}
is transcendental. In fact, when $z:=e^{\,\omega}$ the numbers $^1\!z, ^2\!z,^3\!z,\dotsc$ are algebraically independent, as are $^2\!z, ^3\!z,^4\!z,\dotsc$ when $z:=\alpha$.}\\

We give a  conditional proof of a slightly stronger statement, assuming another conjecture.

The following major open problem in transcendental number theory was stated in the 1960s in a course at Yale given by Lang~\cite[pp. 30--31]{lang}.\\

\noindent{\bf Schanuel's Conjecture (SC).}
\textit{If $\alpha_1,\dotsc,\alpha_n\in\C$ are linearly independent over $\Q$, then the set $\{\alpha_1,\dotsc,\alpha_n,e^{\alpha_1},\dotsc,e^{\alpha_n}\}$ contains at least $n$ algebraically independent numbers.}\\

SC is equivalent to the statement that {\it if $\alpha_1,\dotsc,\alpha_n\in\C$ are $\Q$-linearly independent, then the transcendence degree \emph{trdeg}$_\Q\Q(\alpha_1,\dotsc,\alpha_n,$ $e^{\alpha_1},\dotsc,e^{\alpha_n})$ is at least~$n$}. (See, e.g., Lang \cite[Chapter VIII]{lang2} for the definition and properties of trdeg.)

Note that {\em the case $n=1$ of SC is true} by the \emph{Hermite-Lindemann Theorem} \cite[pp.~6--8]{baker}, which states that $e^\alpha$ is transcendental if $\alpha\neq0$ is algebraic.

For many consequences of SC, see Marques and Sondow~\cite{MS} and its references.
For connections of SC to model theory and exponential algebra, including citations of papers by Ax, Kirby, Macintyre, Wilkie, Zilber, and others, see Waldschmidt~\cite{waldSC} and Wikipedia~\cite{wiki}. 

In the present paper, we introduce a conjecture which is a priori weaker than SC. First, we say that a set $X\subset \C$ is \emph{ $\overline{\Q}$-dependent on} a set $Y\subset \C$ if $\overline{\Q}(X) \subset \overline{\Q}(Y).$\\

\noindent{\bf Schanuel Subset Conjecture (SSC).}
{\it Given $\alpha_1,\dotsc,\alpha_n\in\C$ linearly independent over~$\Q,$ if the set $\{\alpha_1,\dotsc,\alpha_n,$ $e^{\alpha_1},\dotsc,e^{\alpha_n}\}$ is $\overline{\Q}$-dependent on a subset $\{\beta_1,\dotsc,\beta_n\},$
then the numbers $\beta_1,\dotsc,\beta_n$ are algebraically independent.}\\

Note that, by the $\overline{\Q}$-dependence and subset conditions, the fields $\Q(\beta_1,\dotsc,\beta_n)$ and
$\Q(\alpha_1,\dotsc,\alpha_n,e^{\alpha_1},\dotsc,e^{\alpha_n})$ have the same algebraic closure. Thus they also have the same transcendence degree. Now, if SC is true, then we have
\begin{center}
    trdeg$_{\Q}\Q(\beta_1,\dots,\beta_n)$ = trdeg$_{\Q}\Q(\alpha_1,\dots,\alpha_n,e^{\alpha_1},\dots,e^{\alpha_n})\ge n,$
\end{center}
and so $\beta_1,\dotsc,\beta_n$ are algebraically independent. Therefore, \emph{SC implies SSC}.

We do not know whether SSC implies SC.\\

\noindent{\bf Question 1.} Are SC and SSC equivalent?\\

For $n=2$, the answer is yes.\\

\noindent{\bf Proposition 1.}
{\em SSC implies SC for $n=2$.}
\vspace{-2.6mm}
\begin{proof}
If SC is false for $n=2$, then there exist $\Q$-linearly independent complex
numbers $a,b$ such that the field extension $\Q(a,b,e^a, e^b) | \Q$ has transcendence degree 1. Hence 
$\{a, b\} \not\subset\QQ$ (otherwise trdeg$_{\Q}\Q(a,b,e^a, e^b) \ge2$, by the {\em Lindemann-Weierstrass Theorem} \cite[Theorem~1.4]{baker}, which states that \textit{if $\alpha_1,\dotsc,\alpha_k$ are distinct algebraic numbers, then $e^{\alpha_1},\dotsc,e^{\alpha_k}$ are $\QQ$-linearly independent}). Thus we can pick $\beta\in\{a, b\}$
such that $\{\beta\}$ is a transcendence basis for $\Q(a,b,e^a, e^b) | \Q$. By the definition of transcendence basis, the field extension $\Q(a,b,e^a, e^b) | \Q(\beta)$ is algebraic, and then
$$\QQ(a,b,e^a, e^b) \subset \QQ(\beta)  \subset \QQ(\beta,e^{\beta}).$$
Hence the set $\{a,b,e^a, e^b\}$ is $\QQ$-dependent on its subset $\{\beta,e^{\beta}\}$. But then, by SSC, the numbers $\beta,e^{\beta}$ are algebraically independent, contradicting the fact that trdeg$_{\Q}\Q(a,b,e^a, e^b) = 1$.
\end{proof}

Here is our first application of SSC.\\

\noindent{\bf Theorem 1.}
{\it Assume the Schanuel Subset Conjecture. Then Gelfond's Power Tower Conjecture is also true. Moreover, under the weaker hypothesis that $\alpha$ is an algebraic number but not a rational integer, the power tower $^m\!\alpha$ of order $m\ge3$
is transcendental and the numbers $\log\alpha, ^3\!\alpha, ^4\!\alpha,^5\!\alpha,\dotsc$ are algebraically independent.}\\

For example, take $\omega=1$ and $\alpha=1/2.$ Then if SSC holds, the numbers $e,e^{\,e},e^{\,e^{\,e}},\dotsc$ are algebraically independent, and so are the numbers
$$\log2,\ 1/2^{\,1/\sqrt{2}},\ 1/2^{\,1/2^{\,1/\sqrt{2}}},\ 1/2^{\,1/2^{\,1/2^{\,1/\sqrt{2}}}},\dotsc.$$

We give two other consequences of SSC. (They may be known applications of SC, but we have not found them in the literature.) The first states in particular that SSC implies the transcendence of $e^e$ and $\pi^{\pi}$.\\

\noindent{\bf Theorem 2.} \label{PeQe}
{\it If the Schanuel Subset Conjecture is true, then for any non-constant polynomials $P(x),Q(x)\in \QQ[x],$ the numbers $P(e)^{Q(e)}$ and $P(\pi)^{Q(\pi)}$ are transcendental.}\\

Our proof can be adapted to show that SSC also implies the transcendence of $P(\log 2)^{Q(\log 2)}$. On the other hand, there do exist transcendental numbers $T$ for which $T^T$ is algebraic---see \cite[Proposition~2.2]{SM}.

In view of the Gelfond-Schneider Theorem, it is natural to ask: \emph{Which transcendental numbers are not algebraic powers of algebraic numbers?} For instance, $e\neq\alpha^{\beta}$ for any $\alpha,\beta\in\QQ$, since otherwise $e^{1/\beta}=\alpha\in\QQ$ would contradict the Hermite-Lindemann Theorem. Our second application of SSC is that also $\pi\neq \alpha^{\beta}$ and $\log2\neq \alpha^{\beta}$. In fact, we prove a more general statement.\\

\noindent{\bf Theorem 3.} \label{notGS}
{\it Assume that the Schanuel Subset Conjecture is true. Let $\alpha$ and $\beta$ be any algebraic numbers, and let $P(x)\in \QQ[x]$ be a non-constant polynomial. Then $(\alpha^{\beta}-P(e))(\alpha^{\beta}-P(\pi))(\alpha^{\beta}-P(\log2))\neq0.$}\\

We give the proofs of Theorems~1, 2, 3 in Sections~\ref{proof of gelfond1}, \ref{SEC: PeQe}, \ref{SEC: notGS}, respectively.

\section{Proof of Theorem 1} \label{proof of gelfond1}

The proof of Theorem 1 is in two parts.

\begin{proof}[Proof for $z:=e^{\,\omega}$]
As $0\neq\omega\in\QQ$, the Hermite-Lindemann theorem implies that $^1\!z = z = e^{\,\omega}$ is transcendental, and so the statement is true for $k=1$. Now, fix $k>0$ and suppose inductively that the numbers $^1\!z,^2\!\!z,\dotsc,^k\!\!z$ are algebraically independent with $z = e^{\,\omega}$. Then $\omega,\omega\cdot^1\!z,\omega\cdot^2\!z,\dotsc,\omega\cdot\!^k\!z$ are $\Q$-linearly independent. As $^{j+1}\!z  =\,  e^{\omega\cdot\,^j\!z}$, SSC applied to the subset
$$\{^1\!z,^2\!\!z,^3\!\!z,\dotsc,^{k+1}\!\!z\} \subset \{\omega,\omega\!\cdot\!^1\!z,\omega\!\cdot\!^2\!z,\dotsc,\omega\!\cdot\!^k\!z,e^\omega,e^{\omega\cdot\mspace{1mu}^1\!z},e^{\omega\cdot\mspace{1mu}^2\!z},\dotsc,e^{\omega\cdot\mspace{1mu}^k\!z}\}$$
yields the algebraic independence of the numbers $^1\!z,^2\!\!z,\dotsc,^{k+1}\!\!z$. This completes the induction and proves the theorem for $z:=e^{\,\omega}$.
\end{proof}

\begin{proof}[Proof for $z:=\alpha$]
Assuming SSC is true, we show that if $\alpha\in\QQ\setminus\Z$ and $m\ge3$, then the numbers $\log\alpha, ^3\!\!\alpha, ^4\!\!\alpha,\dotsc,^m\!\!\alpha$ are algebraically independent; the proof is in two cases. (A similar proof, but with a single case, shows that if $\alpha\in\QQ\setminus\Q$ and $m\ge2$, then the numbers $^2\!\alpha, ^3\!\!\alpha, \dotsc,^m\!\!\alpha$ are algebraically independent. Details are omitted.)

\noindent
{\it Case 1:} $\alpha\in \Q\setminus \Z$.
Lemma~2.1 in \cite{SM} states that $\alpha^{\alpha}$ is irrational. Thus $1$ and $\alpha^{\alpha}$ are $\Q$-linearly independent, and then so are $\log \alpha$ and $\alpha^{\alpha}\log \alpha$. Since $\alpha$ and $\alpha^{\alpha}$ are algebraic, SSC applied to the subset
$$\{\log \alpha,\alpha^{\alpha^{\alpha}}\} \subset \{\log \alpha,\alpha^{\alpha}\log \alpha, \alpha,\alpha^{\alpha^{\alpha}}\}$$
yields the algebraic independence of $\log \alpha$ and $\alpha^{\alpha^{\alpha}}=\,^3\!\alpha$.

Now, fix $m>3$ and assume inductively that
$\log\alpha, ^3\!\!\alpha, ^4\!\!\alpha,\dotsc,^{m-1}\!\!\alpha$
are algebraically independent. Then in any $\Q$-linear relation
$$\displaystyle\sum_{j=1}^{m-1}a_{j}\cdot\, ^j\!\alpha=0$$
we must have $a_3=\dotsb = a_{m-1}=0$. Since $^1\!\alpha=\alpha\in \Q\setminus \Z$ and $^2\!\alpha=\alpha^{\alpha}\notin \Q$, we also have $a_1=a_2=0$. That implies the $\Q$-linear independence of
$$\{^1\!\alpha\log \alpha,\dotsc,^{m-1}\!\!\alpha\log\alpha\}=\{\log(^2\!\alpha),\dotsc,\log(^m\!\alpha)\}.$$
Then SSC yields the algebraic independence of the subset
\begin{align*}
\{\log(^2\!\alpha),^3\!\!\alpha,^4\!\!\alpha,\dotsc,^m\!\!\alpha\}  \subset \{\log(^2\!\alpha),\dotsc,\log(^m\!\alpha), ^2\!\!\alpha,\dotsc,^m\!\!\alpha\}
\end{align*}
and, hence, since $\log(^2\alpha)=\alpha\log\alpha$, also that of the set $\{\log\alpha,^3\!\!\alpha,^4\!\!\alpha,\dotsc,^m\!\!\alpha\}$. This completes the induction.\\
\noindent
{\it Case 2:} $\alpha\in \QQ\setminus\Q$.
By the Gelfond-Schneider Theorem, $\alpha^\alpha$ is transcendental. Hence $1,\alpha,\alpha^\alpha$ are $\Q$-linearly independent, and then so are $\log \alpha, \alpha\log \alpha, \alpha^\alpha\log \alpha$. Since $\{\alpha\log \alpha, \alpha^\alpha\log \alpha\}\subset \QQ(\log \alpha,\alpha,\alpha^\alpha)$ and $\alpha$ is algebraic, SSC applied to the subset
\begin{center}
$\{\log \alpha, \alpha^\alpha, \alpha^{\alpha^\alpha}\} \subset \{\log \alpha, \alpha\log \alpha, \alpha^\alpha\log \alpha, \alpha, \alpha^\alpha, \alpha^{\alpha^\alpha}\}$
\end{center}
yields the algebraic independence of $\{\log \alpha, \alpha^\alpha, \alpha^{\alpha^\alpha}\}=\{\log \alpha, ^2\!\!\alpha, ^3\!\!\alpha\}$.

Suppose inductively that $\log \alpha, ^2\!\!\alpha,^3\!\!\alpha,\dotsc,^{m-1}\!\!\alpha$ are algebraically independent, where $m>3$. Then any $\Q$-linear relation
$$a_0+\displaystyle\sum_{j=1}^{m-1}a_{j}\cdot\,^j\!\alpha=0$$
implies $a_2=\dotsb =a_{m-1}=0$. Since $^1\!\alpha=\alpha\notin \Q$, we also get $a_0=a_1=0$. That implies the $\Q$-linear independence of
$$\{\log \alpha,^1\!\!\alpha\log \alpha,\dotsc,^{m-1}\!\!\alpha\log\alpha\}=\{\log \alpha,\log(^2\!\alpha),\dotsc,\log(^m\!\alpha)\}.$$
Since
$$\{^1\!\alpha\log \alpha,\dotsc,^{m-1}\!\!\alpha\log \alpha\}\subset \QQ(\log \alpha,^1\!\!\alpha,^2\!\!\alpha,\dotsc,^m\!\!\alpha),$$
we may apply SSC to the subset
\begin{align*}
\{\log \alpha,^2\!\!\alpha,^3\!\!\alpha,\dotsc,^m\!\!\alpha\} \subset \{\log(^1\!\alpha),\dotsc,\log(^m\!\alpha),^1\!\!\alpha,\dotsc,^m\!\!\alpha\}
\end{align*}
and conclude that $\log \alpha,^2\!\!\alpha,^3\!\!\alpha,\dotsc,^m\!\!\alpha$ are algebraically independent.

Thus, in both Cases~1 and~2, the numbers $\log\alpha, ^3\!\alpha, ^4\!\alpha,\dotsc,^m\!\alpha$ are algebraically independent, as desired.
\end{proof}

Results on the arithmetic nature of power towers of $x$ of \emph{infinite} order
$$^\infty x:=\lim_{k\to\infty}\,\!^kx=x^{\,x^{\,x^{\,\cdot^{\cdot^{\cdot}}}}} \qquad\left(e^{-e}\le x\le e^{1/e}\right)$$
can be found in \cite[Appendix]{SM}.

\section{Proof of Theorem~2} \label{SEC: PeQe}

Fix non-constant polynomials $P(x),Q(x) \in \QQ[x]$.

\begin{proof}[Proof that $P(e)^{Q(e)}$ is transcendental] Firstly, let us consider the case $P(x)=x^n$, where $n\ge1$. Since $Q(e)$ is transcendental, $1$ and $nQ(e)$ are $\Q$-linearly independent. Applying SSC to the subset
$$\{e,e^{nQ(e)}\} \subset \{1,nQ(e),e,e^{nQ(e)}\},$$
it follows that $e^{nQ(e)}=P(e)^{Q(e)}$ is transcendental, as claimed.

Now, assume $P(x)\neq x^n$ for any $n\ge1$. We show that $1$ and $\log P(e)$ are $\Q$-linearly independent. Given a $\Q$-linear relation
$a+b\log P(e) = 0$,
by clearing the denominators if necessary, we may assume that $a,b \in \Z$, with $b\geq 0$. Now, $P(e)^be^a-1=0$. If $a\ge0$, then, since $P(x)\neq x^n$ for any~$n\ge0$ and $e$~is not algebraic, the polynomial $P(x)^bx^{a}-1$ must be identically zero; hence $a=b=0$. If $a<0$, then $P(x)^b-x^{-a}=x^{-a}\left(P(x)^bx^{a}-1\right)$ must be the zero polynomial, and again $a=b=0$.

Now, SSC applied to the subset
$$\{e,\log P(e)\} \subset \{1,\log P(e), e,P(e)\}$$
implies that $e$ and $\log P(e)$ are algebraically independent. Hence so are $Q(e)$ and $\log P(e)$. Therefore, the three numbers $1,\log P(e)$, and $Q(e)\log P(e)$ are $\Q$-linearly independent. Applying SSC to the subset
$$\{e,\log P(e),P(e)^{Q(e)}\} \subset \{1,\log P(e),Q(e)\log P(e),e,P(e),P(e)^{Q(e)}\},$$
we get that $P(e)^{Q(e)}$ is transcendental.
\end{proof}

The proof for $P(\pi)^{Q(\pi)}$ will not require two cases, since $1$ and $\log P(\pi)$ are $\Q$-linearly independent even when $P(x)=x^n$.

\begin{proof}[Proof that $P(\pi)^{Q(\pi)}$ is transcendental] 
Note first that $i\pi$ and $\log P(\pi)$ are $\Q$-linearly independent, for if there existed a $\Z$-relation
   $ai\pi+ b\log P(\pi)=0$ with $b>0$, then $P(\pi)^b=(-1)^a$ would be algebraic, contradicting the transcendence of~$\pi$.

Applying SSC to the subset
$$\{i\pi,\log P(\pi)\} \subset \{i\pi,\log P(\pi),-1,P(\pi)\},$$
we get that $i\pi$ and $\log P(\pi)$ are algebraically independent. Then the set
$$\{i\pi, \log P(\pi), Q(\pi)\log P(\pi)\}$$ is $\Q$-linearly independent, and SSC applied to this subset of
$$\{i\pi, \log P(\pi), Q(\pi)\log P(\pi),-1,P(\pi),P(\pi)^{Q(\pi)}\}$$
yields the desired result.
\end{proof}

\section{Proof of Theorem 3} \label{SEC: notGS}

It suffices to show that if SSC is true and $\omega\neq0$ is algebraic, then the numbers $P(e)^{\omega}$, $P(\pi)^{\omega}$, and $P(\log2)^{\omega}$ are all transcendental. Before giving the proof, we first prove an unconditional lemma.\\

\noindent{\bf Lemma 1.}
{\it If $Q(x)\in \QQ[x]\setminus \{0,x^n:n=0,1,2,\dotsc\}$, then $\log Q(e)$ is transcendental.}

\begin{proof} Suppose on the contrary that $\alpha:=\log Q(e)\in\QQ$. If $Q(x)=\sum_{k=0}^n a_kx^k$\\ then we have the relation
$$
a_0+a_1e+\dotsb +a_ne^n-e^{\alpha}=0.
$$
The Lindemann-Weierstrass Theorem implies first that $\alpha\in \{0,1,\dotsc,n\}$, and then that $a_\alpha=1$ and $a_k=0$, for $k\neq\alpha$. But then $Q(x)=x^\alpha$, contradicting the hypothesis. Therefore $\log Q(e)\not\in\QQ$.
\end{proof}

Since the hypotheses of Theorem \ref{notGS} imply that $P(e),P(\pi),$ and $P(\log2)$ are transcendental, in order to prove that if $0\neq\omega\in \QQ$, then $P(e)^\omega,P(\pi)^\omega,$ and $P(\log2)^\omega$ are also transcendental, we only need to consider the case $\omega\in \QQ\setminus\Q$.

\begin{proof}[Proof that $P(e)^\omega\not\in\QQ$]
Since $e^{n\omega}$ is transcendental for $n=1,2,\dotsc$, we may assume $P(x)\neq x^n$. Then Lemma~1 implies $\log P(e)$ is transcendental, which in turn implies the $\Q$-linear independence of $1,\log P(e),\omega\log P(e)$. Now, by SSC applied to the subset
$$
\{e,\log P(e),P(e)^{\omega}\}\subset \{1,\log P(e),\omega\log P(e),e, P(e),P(e)^{\omega}\},
$$
we get, in particular, the transcendence of $P(e)^{\omega}$.
\end{proof}

\begin{proof}[Proof that $P(\pi)^\omega\not\in\QQ$.]
In the second part of the proof of Theorem~2, we proved that $i\pi$ and $\log P(\pi)$ are algebraically independent. Since $\omega$ is irrational, the set $\{i\pi,\log P(\pi), \omega\log P(\pi)\}$ is $\Q$-linearly independent. Now, we can apply SSC to the subset
$$
\{\pi,\log P(\pi),P(\pi)^{\omega}\}\subset \{i\pi,\log P(\pi), \omega\log P(\pi),-1,P(\pi),P(\pi)^{\omega}\}
$$
and conclude that $P(\pi)^{\omega}$ is transcendental.
\end{proof}

\begin{proof}[Proof that $P(\log2)^\omega\not\in\QQ$.]
The numbers $\log 2$ and $\log P(\log 2)$ are $\Q$-linearly independent. In fact, any $\Q$-relation $a\log 2+b\log P(\log 2)=0$ implies that $P(\log 2)^b=2^{-a},$ and then $a=b=0$ by the transcendence of $P(\log 2)$. By SSC applied to the subset
$$
\{\log 2,\log P(\log 2)\}\subset \{\log 2,\log P(\log 2),2,P(\log 2)\},
$$
we have that $\log 2,\log P(\log 2)$ are actually algebraically independent, and so the set $\{\log 2, \log P(\log 2),\omega\log P(\log 2)\}$ is $\Q$-linearly independent. Again by SSC, applied to the subset
\begin{align*}
&\{\log 2, \log P(\log 2), P(\log 2)^{\omega}\}\\
\subset &\{\log 2,\log P(\log 2),\omega\log P(\log 2),2, P(\log 2), P(\log 2)^{\omega}\},
\end{align*}
we get that $\log 2, \log P(\log 2), P(\log 2)^{\omega}$ are algebraically independent. The theorem follows.
\end{proof}

We leave it as an exercise to use similar arguments to show, under the assumption of SSC, that each of the numbers $e^e$ and $e+\pi$ is also not an algebraic power of an algebraic number.

\bibliographystyle{amsplain}

\end{document}